\documentclass[11pt]{amsart}
\usepackage{setspace, amsmath, amsthm, amssymb, amsfonts, amscd, epic, graphicx, ulem, dsfont}
\usepackage[T1]{fontenc}
\usepackage{multirow}
\usepackage{bbm}
\usepackage{enumerate}

\makeatletter \@namedef{subjclassname@2010}{
  \textup{2010} Mathematics Subject Classification}
\makeatother

\newtheorem{thm}{Theorem}[section]

\newtheorem{lem}[thm]{Lemma}

\theoremstyle{remark}
\newtheorem*{rema}{Remark}

\newtheorem{exa}[thm]{\textbf{Example}}

\theoremstyle{definition}

\begin{document}

\title{An Improvement of Reid Inequality}
\author[M. H. Mortad]{Mohammed Hichem Mortad}

\dedicatory{}
\thanks{}
\date{}


\address{Department of
Mathematics, University of Oran 1, Ahmed Ben Bella, B.P. 1524, El
Menouar, Oran 31000, Algeria.\newline {\bf Mailing address}:
\newline Pr Mohammed Hichem Mortad \newline BP 7085 Seddikia Oran
\newline 31013 \newline Algeria}

\email{mhmortad@gmail.com, mortad@univ-oran.dz.}

\begin{abstract}In this short note, we improve the famous Reid Inequality related
to linear operators.
\end{abstract}

\maketitle

\section{Main Result}
First, assume that readers are familiar with notions and result on
$B(H)$. We do recall a few definitions and results though:
\begin{enumerate}
  \item Let $A\in B(H)$. We say that $A$ is positive (we then write $A\geq 0$) if
  \[<Ax,x>\geq0,~\forall x\in H.\]
  \item For every positive operator $A\in B(H)$, there is a unique positive $B\in B(H)$
  such that $B^2=A$. We call $B$ the positive square root of $A$.
  \item The absolute value of $A\in B(H)$ is defined to be the
  (unique) positive square root of the positive operator $A^*A$.
  We denote it by $|A|$.
  \item We recall that $A\in B(H)$ is called hyponormal if $AA^*\leq A^*A$.
\end{enumerate}

The inequality of Reid which first appeared in
\cite{Reid-Inequality} is recalled next:

\begin{thm}\label{Redi classic THM}
Let $A,K\in B(H)$ such that $A$ is positive and $AK$ is
self-adjoint. Then
\[|<AKx,x>|\leq ||K||<Ax,x>\]
for all $x\in H$.
\end{thm}

\begin{rema}
As shown in e.g. \cite{Lin-Reid-Furuta-inequality}, Reid Inequality
is equivalent to the operator monotony of the positive square root
on the set of positive operators.
\end{rema}

Many generalizations of Theorem \ref{Redi classic THM} are known in
the literature from which we only cite \cite{Kitanneh-REid-et al}
and \cite{Lin-Reid-Furuta-inequality}.

In an earlier version of this paper (see \cite{Mortad-Reid-Normal!
FIRST VERSION}), the author showed the following:

\begin{thm}\label{Main THM}
Let $A,K\in B(H)$ such that $A$ is positive and $AK$ is normal. Then
\[|<AKx,x>|\leq ||K||<Ax,x>\]
for all $x\in H$.
\end{thm}

Can we go to $AK$ being hyponormal? The answer is no as seen next:

\begin{exa}
Let $S$ be the shift operator on $\ell^2$. Setting $A=SS^*$, we see
that $A\geq0$. Now, take $K=S$ (and so $\|K\|=1$). It is clear that
$AK=SS^*S=S$ is hyponormal. If Reid Inequality held, then we would
have
\[|<Sx,x>|\leq <SS^*x,x>=\|S^*x\|^2\]
for each $x\in\ell^2$. This inequality clearly fails to hold for all
$x$. Indeed, taking $x=(2,1,0,0,\cdots)$, we see that
\[|<Sx,x>|=2\leq \|S^*x\|^2=1\]
which is impossible.
\end{exa}

The good news is that Reid Inequality can yet be improved as it
holds if $AK$ is co-hyponormal, that is, if $(AK)^*$ is hyponormal.
This comes after a discussion with a fellow student (Mr S. Dehimi):

\begin{thm}\label{Main THM}
Let $A,K\in B(H)$ such that $A$ is positive and $(AK)^*$ is
hyponormal. Then
\[|<AKx,x>|\leq ||K||<Ax,x>\]
for all $x\in H$.
\end{thm}

The proof relies on the following result:

\begin{lem}\label{Kittaneh Lemma abos valu hypo}(\cite{Kitanneh-REid-et
al}) Let $A\in B(H)$ be hyponormal. Then
\[|<Ax,x>|\leq <|A|x,x>.\]
\end{lem}

Now, we give the proof of Theorem \ref{Main THM}.

\begin{proof}
  The inequality is evident when $K=0$. So, assume that
$K\neq0$. It is then clear that $\frac{K}{\|K\|}$ satisfies
\[KK^*\leq \|K\|^2I.\]
Hence
\[|(AK)^*|^2=AKK^*A\leq \|K\|^2A^2\]
or simply $|(AK)^*|\leq \|K\|A$ after passing to square roots.

Now, for all $x\in H$
\[|<AKx,x>|=|<x,(AK)^*x>|=|\overline{<(AK)^*x,x>}|=|<(AK)^*x,x>|.\]

Since $(AK)^*$ is hyponormal, Lemma \ref{Kittaneh Lemma abos valu
hypo} combined with $|(AK)^*|\leq \|K\|A$ give
\[|<AKx,x>|=|<(AK)^*x,x>|\leq |<|(AK)^*|x,x>|\leq \|K\|<Ax,x>\]
and this marks the end of the proof.
\end{proof}

\end{document}